\documentclass[dvips]{imsart}
\usepackage[OT1]{fontenc}
\usepackage{amsthm,amsmath}
\usepackage{amsfonts}
\usepackage{amsbsy}
\usepackage{subfigure}
\usepackage{wrapfig}
\usepackage{longtable}
\usepackage{amssymb}
\usepackage{latexsym}
\usepackage{euscript}
\usepackage{fancyhdr}
\usepackage{epsfig}
\usepackage{graphicx}
\usepackage{subfigure}
\usepackage{fancybox}
\usepackage{varioref}
\usepackage{enumerate}
\usepackage{verbatim}
\usepackage{makeidx}
\usepackage[english]{babel}
\usepackage[applemac]{inputenc}

\newcommand{\R}{\mathbb{R}}
\newcommand{\E}{\mathbb{E}}

\newcommand{\1}{{\sf \hspace*{0.9ex}}\rule{0.15ex}{1.6ex}\hspace*{-1ex} 1}

\newcommand{\Argmax}{Argmax}

\newcommand{\X}{\mathcal{X}}

\arxiv{math.PR/00000000}

\startlocaldefs
\theoremstyle{plain}
\newtheorem{Lemme}{Lemma}[section]
\newtheorem{proposition}{Proposition}[section]
 
 \newtheorem{theoreme}{Theorem}[section]

\endlocaldefs

\begin{document}

\begin{frontmatter}
\title{Fast rate of convergence in high dimensional linear discriminant analysis}
\runtitle{High dimensional discriminant analysis}

\begin{aug}
\author{ \snm{Girard}\fnms{R.}\ead[label=e1]{robin.girard@mines-paristech.fr}}

\address{Mines Paristech, Sophia Antipolis, France\\
\printead{e1}}

\runauthor{R. Girard}
\affiliation{Mines Paritech}
\end{aug}
\begin{abstract}
This paper gives a theoretical analysis of high dimensional linear discrimination of Gaussian data. We study the excess risk of linear discriminant rules. We emphasis on the poor performances of standard procedures in the case when dimension $p$ is larger than sample size $n$. The corresponding theoretical results are non asymptotic lower bounds. On the other hand, we propose two discrimination procedures based on dimensionality reduction and provide associated rates of convergence which can be $O\left( \frac{\log(p)}{n} \right )$ under sparsity assumptions. Finally all our results rely on a theorem that provides simple sharp relations between the excess risk and an estimation error associated to the geometric parameters defining the used discrimination rule.
 \end{abstract}



\end{frontmatter}

\section{Introduction}

 In the binary classification problem, the aim is to recover the unknown class $y\in \{0,1\}$ associated to an observation $x\in \X= \R^p$. In other words, we seek a classification rule, also called classifier: a measurable $g:\X \rightarrow \{0,1\}$. This rule gives a wrong classification for the observation $x\in \R^p$ if $g(x)\neq y$. The underlying probabilistic model, that allows us to measure the performances of a classification rule $g$, is set by a distribution $P$ on $\X\times \{0,1\}$ with conditional probability $P_k()=P(.\times \{k\})$ ($k=0,1$). In this framework, under a uniform prior, the probability of misclassification is defined by 
\[\mathcal{C}(g)=\frac{1}{2}\left ( P_1(X\notin g^{-1}({1}))+P_0(X\notin g^{-1}({0}))\right ).\]

 In this paper we consider the case when $P_0$ and $P_1$ are gaussian with mean $\mu_0$ and $\mu_1$ respectively and with common covariance $C$. Since then, when $\X=\R^p$, the Bayes rule, i.e the classification rule $g^*$ that minimizes $\mathcal{C}(g)$, is given by 
\begin{equation}\label{bayes}
g^{*}(x)=\left\{\begin{array}{cc}1 & \text{if }\langle F_{10},x-s_{10}\rangle_{\R^p} \geq 0\\ 0 & \text{ otherwise }\end{array}\right .
\end{equation}
\[\text{where }F_{10}=C^{-}(\mu_1-\mu_0),\; s_{10}=\frac{\mu_1+\mu_0}{2},\]
 $C^-$ is the generalized inverse\footnote{If $C$ is a semi-positive definite matrix, one can define the associated generalized inverse, also called Moore-Penrose pseudo-inverse : $C^-$. This generalised inverse $C^-$ arises from the decomposition $\R^p=Ker(C) \oplus Ker(C)^{\bot}$. On $Ker(C)$, $C^-$ is null, ad on $Ker(C)^{\bot}$, $C^-$ equals the inverse of $\tilde{C}=C_{|Ker(C)^{\bot}}$ (  i.e $\tilde{C}$ is the restriction of  $C$ to $Ker(C)^{\bot}$).}  of $C$ and $\langle\;,\; \rangle_{\R^p}$ is the euclidian inner product of $\R^p$. Since $\mu_1,\mu_0$ and $C$ are unknown, $g^*$ is unknown. Assume that one observes two independent samples $X^0=(X_1^0,\dots,X_{n_0}^0)$ $X^1=(X_1^1,\dots,X_{n_1}^1)$ of $\X$ valued i.i.d observations with probability distribution $P_0$ or $P_1$, respectively. One can use empirical rules $\hat{g}_{n_{0},n_1}$ based on the observations $X^0,X^1$ to mimic $g^*$. When one assumes that $P_1$ and $P_0$ are gaussian with the same covariance, it becomes natural to search for a classification rule $g:\R^p\rightarrow \{ 0,1\}$  given by
\begin{equation}\label{gdef}
g(x)=\left\{\begin{array}{cc}1 & \text{if }\langle \hat{F}_{10},x-\hat{s}_{10}\rangle_{\R^p} \geq 0\\ 0 & \text{ otherwise }\end{array}\right.
\end{equation} 
where $\hat{F}_{10},\hat{s}_{10}\in \R^p$ have to be estimated from the observations $X^0,X^1$.\\
 
   A standard way of assessing the quality of a decision rule $\hat{g}_{n}$ (where $n=n_{1}+n_0$)  is to give an upper bound on $\E[\mathcal{C}(\hat{g}_{n})-\mathcal{C}(g^*)]$. A classification rule $\hat{g}_{n}$ is said to be consistent if this last quantity converges to zero when $n\rightarrow \infty$. In this paper,  we are interested in the case where $p>>n=$ ($p$ is the dimension of $\X$), and our aim is twofold. First, we give two procedures to achieve the fast rate of convergence. These procedures rely on a dimensionality reduction. Second, we give lower bounds on the excess risk to show that standard procedures (such as the Fisher discriminant analysis) fail in high dimension (when $p>n$).  These lower bounds are given as a function of the sample size $n$ and the dimension $p$. They are not asymptotic lower bounds since these bounds remain valid  for all the cases when $p>n$. \\
 
   Let us introduce some notations that will be used throughout this paper. If $P$ is a probability measure on $\R^p$ with finite second order moment and $u,v\in \R^p$, $\|v\|_{L_2(P)}$ will stand for the $L_2(P)$ norm\footnote{Let us recall that the $L_2(P)$ norm of $f:x\in \R^p\rightarrow f(x)$ is defined by $\|f\|^{2}_{L_2(P)}=\int_{x\in \R^p}f^2(x)P(dx)$} of $x\in \R^p\rightarrow \langle v , x \rangle_{ \R^p}$, and $\langle u , v\rangle_{L_2(P)}$ will stand for the associated scalar product. This scalar product induces a geometry in $\R^p$, the associated angle in $L_{2}(P)$ between $u$ and $v$ will be denoted by $\alpha_{L_{2}(P)}(u,v)$. In the rest of the paper, $P_C$ will stand for a gaussian centered measure with covariance $C$.

Our main result in this paper is Theorem \ref{thbp}. There, we see that when
 \begin{enumerate}
 \item $F_{10}$ has a finite number of non null components (sparsity assumption)
 \item $\|F_{10}\|_{L_2(P_C)}$ is lower bounded ( strict margin assumption)
 \end{enumerate}
 then the procedure we are proposing achieves the rate $\log(p)/n$. Finally, our theorem also shows identical rate of convergence for other types of sparsity assumption and margin assumptions.\\

There is a large body of literature about lower bounds on the excess risk in the classification framework, one can see for example \cite{Massart:2006qq,Audibert:2006fk,Mammen:1998uq,Yang:1998qf,optaggr}. These articles are mainly dedicated to the problem of finding the minimax rate of convergence in certain classes of classifiers.  These classes cannot be adapted to our case.  Moreover, we do not search minimax lower bounds.\\

 The classification rule we propose is a linear discriminant analysis with a dimensionality reduction procedure. This type of discrimination  procedure in a high dimensional gaussian framework has been investigated in \cite{donoho_feature_2009,Hall:2010fk,Tibshirani:2003im,bickel:2004fk,Fans:2007vn} and our work is in line with these papers.  The main improvement we give is that the full proposed procedures (including the use of a data dependent threshold)  come with a rate of convergence that can be the fast rate under a wide range of sparsity assumptions. In our work we relate classification error and error made while estimating $F_{10}$ and $s_{10}$, also our work is related to the area of plugin classification. Our theoretical development is centered on Theorem \ref{Th1}. There, we give a bound exhibiting a good relation (sharp lower and upper bound) between the estimation error of $F_{10}$ and the excess risk and this has never been investigated. 

This paper is structured as follows. In Section 2 we give finite sample lower bounds showing how bad are standard procedure for finding $F_{10}$ when $p>>n$. In section 3 we give two algorithms to overcome these problems together with associated theoretical results and numerical  experiments. The proofs, and the statement of Theorem \ref{Th1} are postponed to the Annex.

 \section{Inconsistency of standard procedure when there are more variables than observations}

Within the learning set, we observe two independent samples $X^0=(X_1^0,\dots,X_{n_0}^0)$, $X^1=(X_1^1,\dots,X_{n_1}^1)$ of $\R^p$ valued i.i.d observations with probability distribution $P_0$ or $P_1$, respectively. The following proposition illustrates the inconsistency of standard procedures when $p>n=n_1+n_0$. 

\begin{proposition}\label{prop:un}
For $k=0,1$, let $\hat{C}^k$ be the empirical covariance matrix of $X^k$, and  $\hat{\mu}_k$ be the empirical mean of $X^k$. Let us define 
\[\hat{C}=\frac{1}{n-1}\left ((n_{0}-1)\hat{C^0}+(n_1-1)\hat{C}^1\right ),\] 
and let $\hat{s}_{10}$ be any estimator of $s_{10}\in \R^p$. 
\begin{itemize}
\item If $\hat{F}_{10}= \hat{C}^-m_{10}$, then, the classification rule $g$ defined by (\ref{gdef}) leads to   
  \[\E_{P^{\otimes n}}[\mathcal{C}(g)-\mathcal{C}^*]\geq \left (1-\sqrt{\frac{n}{p}}\right ) \frac{\|F_{10}\|_{L_2(P_C)}}{2\sqrt{2\pi}}e^{-\frac{5\|F_{10}\|^2_{L_2(P_C)}}{8}}.\]
\item If $\hat{F}_{10}=C^{-}(\hat{\mu}_1-\hat{\mu}_0)$, then, the classification rule $g$ defined by (\ref{gdef}) leads to
  \[\E_{P^{\otimes n}}[\mathcal{C}(g)-\mathcal{C}^*]\geq \left (1-\frac{\sqrt{n}\|F_{10}\|_{L_2(P_C)}+1}{\sqrt{p-2}}\right )\frac{\|F_{10}\|_{L_2(P_C)}}{2\sqrt{2\pi}} e^{-\frac{5\|F_{10}\|^2_{L_2(P_C)}}{8}}\] 
\end{itemize}
  
\end{proposition}



 \paragraph{General comments}
First, we note that $d=\frac{\|F_{10}\|_{L_2(P_C)}}{2}$ is related to the $L_1$ distance bewteen $P_0$ and $P_1$ through this known equality:
\[d_1(P_1,P_0)=\int |dP_1-dP_0|=\Phi(-d)-\Phi(d).\]
where $\Phi(x)$ is the cumulative distribution function of a real gaussian random variable with mean zero and variance one. Hence $d_1(P_1,P_0)\sim d$ when $d$ tends to zero. In this case, the preceding lower bound is tight since $\mathcal{C}(g)-\mathcal{C}^*\leq d_1(P_1,P_0)$.
When $d_1(P_1,P_0) \rightarrow 1$, $d\rightarrow \infty$ and  
 \[d_1(P_1,P_0)\sim 1-\frac{e^{-\frac{d^2}{2}}}{d\sqrt{2\pi}}.\]

As a particular application of this proposition, we see that the Fisher Rule is not consistent when $p>>n$, which was already given in \cite{bickel:2004fk}. However, our result is stronger: we can even say that if there exists $1>c>0$ such that $\frac{n}{p}<c$, then the Fisher rule is not consistent.\\

\paragraph{Structural assumption.}
The preceding proposition suggests that in the problem of estimating $F_{10}$ to construct a consistent rule $g$ (as given by Equation \ref{gdef}), when $p>>n$,  a structural assumption on $(C^{-})^{1/2}(\mu_1-\mu_0)$ has to be made (by abuse of notation we will write $C^{-1/2}m_{10}$ in the remaining of the paper).
Indeed, from point $2$ of the proposition, if there exists $0<r<R$ such that $R>\|F_{10}\|^2_{L_2(P_C)}\geq r$, then, uniformly on all the possible values of $\mu_1$ and $\mu_0$, the excess risk can converge to zero only if  $\frac{p}{n}$ tends to $0$. Recall that if no a priori assumption is done on $m_{10}$, $\bar{m}_{10}$ is the best estimator of $m_{10}$ with respect to the quadratic loss: $\bar{m}_{10}=Argmin_{f(X^1,X^0)}\E[\|m_{10}-f(X^1,X^0)\|_{\R^p}^2]$. \\

In the literature of high dimensional classification, the mean difference vector $m_{10}=\mu_1-\mu_0$ is commonly believed to be sparse (see \cite{Fans:2007vn}). In this paper $C^{-1/2}m_{10}$ is assumed to be sparse. Intuitively, the sparsity assumption permits to bound the dimension of $\R^p$ subspace for which the classification can be performed efficiently, and it is sufficient but not necessary to relate this space to the sparsity of $m_{10}$ only. Indeed, there can be a direction $e\in \R^p$ such that $e =argmax_{\|e\|=1}\langle m_{10},e_i\rangle$ but $e =argmin_{\|e\|=1}\langle C^{-1/2}m_{10},e_i\rangle$ and it is natural to take into account the overall dispersion of the data as well as the mean difference vector.\\

 Theoretically, the choice of a sparsity assumption on $C^{-1/2}m_{10}$ is enlighten by Theorem \ref{Th1}. Indeed, this Theorem implies that if $C_-<\|F_{10}\|_{L_2(P_C)}<C_+$ (for $C_-,C_+>0$), there exists $0<C_1<C_2$ such that 
\[C_1 \alpha^2\leq \mathcal{C}(g)-\mathcal{C}^*\leq C_2 \alpha^2+f(\hat{s}_{10})\] 
where $\alpha=\alpha_{L_2(P_C)}(F_{10},\hat{F}_{10})$ is the angle between $F_{10}$ and $\hat{F}_{10}$ in the geometry of $L_2(P_C)$ and $f:\R^p\rightarrow \R_+$ with  $f(s_{10})=0$. This explains why an assumption on the sparsity of $F_{10}$ in $L_2(P_C)$ (or a sparsity of $C^{-1/2}m_{10}$) is more suitable. \\

 The structural assumption on $C^{-1/2}m_{10}$ can be a consequence  of structural assumptions on $\mu_1-\mu_0$ and on $C$. Many works, based on model selection or aggregation have already been done to define proper structural assumption for the estimation of $C$, see for example \cite{Bickel:2007fk} and the reference therein. Those works are dedicated to the problem of estimating $C$ with a Hilbert-Schmidt error measure, and yet do not give results in the classification framework. In addition, we will see in next section that it is not necessary to estimate all the parameters of $C$ but that one only need to estimate $F_{10}$ which has only $p$ parameters. \\
If a structural assumption is done on $C$, it has to be linked with a statistical assumption. For example reducing the number of parameters to estimate can be done with a  stationarity (or quasi stationarity, as in \cite{bestbasis})  assumption. If $C$ is Toeplitz (i.e $C_{ij}=c(i-j)$ with $c:\mathbb{Z}\rightarrow \R$ a $p$-perioric sequence) it is a circular convolution operator which is known to be diagonal in the discrete Fourier basis $(g^m)_{0\leq m<p}$ defined by: 
\[(g^{m})_k=\frac{1}{\sqrt{p}}\exp\left (\frac{2i\pi mk}{p}\right ).\]   
 This is a generalization (to the infinite dimensional framework) of this harmonic analysis result that is used in Bickel et Levina \cite{bickel:2004fk} and combined with approximation in \cite{bestbasis}. Using this type of assumption, the covariance matrix can be searched in the set of diagonal matrices. Let us note that the use of harmonic analysis and stationarity in curve classification can become a wide field of interest as soon as one considers the larger class of group stationnary-processes (see \cite{Yazici:2004mz}) or semi-group stationnary processes (see \cite{Girardin:2003rt}). \\

However, we believe that making directly a structural assumption on $C^{-1/2}m_{10}$ is more suitable in the case or our classification problem. In the estimation of a high dimensional vector problem, finding suitable structural assumption has been studied extensively (see for example \cite{Candes:2006lk}). In this paper, we limit our work to $l^q$ bodies for $0<q<2$. Let $P_C$ be a gaussian measure on $\R^p$ with full rank covariance, for $0<q<2$ let us define $l^q(R,P_C)$ the $l^q$ ball of $L_{2}(P_C)$ with radius $R>0$ by  
\[
l^q(R,P_C)=\left \{v\in \R^p\; :\;  \|C^{1/2}v\|^q_q\leq R^q \right \},
\]
where $\|x\|^q_q=\sum_{i=1}^p|x[i]|^q$ for any $x\in \R^p$. 
 For a well chosen orthonormal basis of $\R^p$, knowing that $F_{10}\in l^{q}(R,P_C)$  for $0<q<2$ will be used (see next Section) to construct a consistent estimator of $F_{10}$.

\section{Fast rate of convergence for linear discrimination rule}

In this section we suppose that $C$ is diagonal, and use the notation $\sigma^2[i]=C[i,i]$. The learning set $(X^k_{j})_{k=0,1,\;j=1,\dots,n_k}$ is separated in two parts, part $A$ and part $B$, with equal size: 
\[\text{ Part }A =(X^k_{j})_{k=0,1,\;1\leq j<n_k/2} \text{ and Part }B =(X^k_{j})_{k=0,1,\;n_k/2\leq j \leq n_k}.\] For $k=0,1$ let $\bar{\mu}^A_{k}$ (resp $\bar{\mu}^B$) be the empirical mean of the learning data from part $A$ (respectively from part $B$) and class $k$. For $i=1,\dots,p$, $k=0,1$, let $\hat{\sigma}_k^2[i]$ be the empirical (unbiased) variance of the $i^{th}$ feature within the learning data from part $B$: $(X^k_{j}[i])_{k\;, n_k/2\leq j \leq n_k}$ and define $\hat{\sigma}^2[i]=\frac{1}{n-1}((n_0-1)\hat{\sigma}_0^2[i]+(n_1-1)\hat{\sigma}_1^2[i])$. Now, let us define 
\begin{equation}\label{s10m10hat}
\hat{s}_{10}=\frac{\bar{\mu}^A_1+\bar{\mu}^A_0}{2},\;\; \bar{m}_{10}=\bar{\mu}^B_1-\bar{\mu}^B_0,\; \hat{\sigma}=(\hat{\sigma}[i])_{i=1\dots,p},\text{ and } \tilde{F}_{10}=(\bar{m}_{10}[i]/\hat{\sigma}^2[i])_{i=1\dots,p}. 
\end{equation}
 We recall that in this paper, $n=n_1+n_0$.
We will note
\begin{equation}
\Omega_q(R,r)=\left \{(P_1,P_0)\in \mathcal{P} \; \text{ s.t } F_{10}\in l^q(R,P_{Cov(P_1)}),\;\;\|F_{10}\|_{L_{2}(P_C)}\geq r\right\}
\end{equation}
where $\mathcal{P}$ is the set of pairs $(P_1,P_2)$ of gaussian probability distribution on $\R^p$ with $cov(P_1)=cov(P_2)$. 
\paragraph{Definition of procedures.}

We propose two discrimination procedures. The first one is simpler and comes with a more complete theoretical result while the second one is more sophisticated but requires further theoretical work. Both use the discrimination procedure $g$ (defined by Equation (\ref{gdef})) with $\hat{s}_{10}$ defined by Equation \ref{s10m10hat}. In both cases $\hat{F}_{10}$ is evaluated upon a dimensionality reduction step

\begin{equation}\label{hatdeF}
\hat{F}_{10}[i]=\left\{\begin{array}{cc} \tilde{F}_{10}[i] & \text{ if } i\in \hat{I}  \\0 & \text{ else }\end{array}\right. 
\end{equation}

Given that we know $\hat{I}$, the preceding rule is a rephrasing of the feature annealed independence rule introduced in \cite{Fans:2007vn} in the case when group $0$ and group $1$ have equal variance. The proposed methods differ by procedure used to construct $\hat{I}$ (even if the result we give in Theorem \ref{Th1} applies to the case when $C$ is not diagonal.). We propose to use two simple procedures borrowed from the thresholding estimation literature (Procedure $1$ and $2$ below) for selecting the subset $\hat{I}$ in order to estimate the normal vector $F_{10}$ to the optimal separating hyperplane. Procedure $3$ is the thresholding procedure proposed in\cite{Fans:2007vn}.  

{\it Procedure 1 : universal dimensionality reduction.} In the first procedure $\hat{I}$ is given by 
\begin{equation}\label{universal}
\hat{I}^U =\left \{i\in \{1,\dots, p\}:\; \left |\frac{\bar{m}_{10}[i]}{\hat{\sigma}[i]}\right |>\sqrt{2\frac{\log(p)}{n}}\right \}
\end{equation}
This can be seen as a thresholding estimation of $C^{-1/2}m_{10}$ with a universal threshold (see for example \cite{Donoho:1998fk}). The next procedure relies on the same idea with a false discovery rate thresholding procedure.\\

{\it Procedure 2 : False discovery rate control dimensionality reduction.} In the second procedure $\hat{I}$ is given by 
\begin{equation}\label{selectFDR}
\hat{I}^{FDR} =\left \{i\in \{1,\dots, p\}:\; \left |\frac{\bar{m}_{10}[i]}{\hat{\sigma}[i]}\right |>\lambda_{FDR}\right \}
\end{equation}
where $\lambda_{FDR}$ is a data dependent threshold chosen with the Benjamini and Hocheberg procedure \cite{BH95} for control of the false discovery rate (FDR) of the following multiple hypothesis :
\begin{equation}\label{multipletest}
\forall i=1,\dots,p \;\;\; H_{0i}\;:\; E\left [\frac{\bar{m}_{10}[i]}{\hat{\sigma}[i]}\right ]=0\;\;\;:\text{ Versus } H_{1i}\;:\; E\left [\frac{\bar{m}_{10}[i]}{\hat{\sigma}[i]}\right ]\neq0
\end{equation}
This procedure is as follows. Let us define $T[i]=\frac{\bar{m}_{10}[i]}{\hat{\sigma}[i]}$. The $(|T[i]|)_i$ are ordered in decreasing order :
\[
|T[(1)]|\geq \dots \geq |T[(p)]|\text{ and }  \lambda_{10}^{FDR}=|T[(k_{10}^{FDR})]|
\]
\[
\text{ where } k_{10}^{FDR}=\max\left \{k \in \{1,\dots,p\}\; : \; |T[(k)]|\geq \sqrt{\frac{1}{n}}z\left (\frac{b_pk}{2p}\right )\right \},
\]
  $z(\alpha)$ is the quantile of order $\alpha$ of a standardized gaussian random variable and $b_p\in [0,1/2[$ is under bounded by $\frac{c_0}{\log p}$ where $c_0$ is a positive constant (which does not depend on $p$).\\

 The procedure can also be seen as a thresholding estimation of $C^{-1/2}m_{10}$, but with a FDR threshold (see \cite{ABDJ05}). There are a lot of thresholding procedures in the literature today and others could be used. The universal threshold is the first that has appeared and the simplest. The FDR threshold is one of the most efficient and adaptive one. In addition, in our problem, it can lead to an interesting statistical rephrasing of the procedure. Indeed, the multiple hypothesis given by Equation (\ref{multipletest}) are connected heuristically with $\forall i=1,\dots,p$ 
\[
H_{0i}\;:\; \text{ the ratio variance inter/variance intra is null in direction i}
\]
Versus
\[H_{1i}\;:\; \text{ the ratio variance inter/variance intra is not null in direction i}.\]
Hence our procedure can be rephrased in to step : 
\begin{enumerate}
\item Make a "vertical analysis of the variance" to select the directions $i\in \hat{I}$ in which the data are well separated (i.e $(C^{-1/2}m_{10})[i]$ is large)  
\item Perform a standard discriminant analysis in the space spanned by the directions chosen in step $1$.
\end{enumerate}

{\it Procedure 3 :  threshold choice from \cite{Fans:2007vn}.} In procedure $3$, $\hat{I}^{FAIR}$ is computed the same way as $\hat{I}^{FDR}$ replacing $k^{FDR}$ by:
\begin{equation}\label{equ:FAIR}
k^{FAIR}=\Argmax_{m=1,\dots,p} \left \{\frac{1}{\max_{i\leq m}\hat{\sigma}^2[(i)]}\frac{n\left (\sum_{i=1}^mT^2_{(i)}+m(1/n_1-1/n_0)\right )^2}{mn_1n_0+n_1n_0\sum_{i=1}^mT^2_{(i)}}\right \}.
\end{equation}
{\it Procedure 4 :  Higher Criticism from \cite{donoho_feature_2009}.} In procedure $4$, $\hat{I}^{HC}$ is computed with the higher Criticism procedure \cite{donoho_feature_2009}: with $\lambda_{HC}=|T[(k^{CH})]$ where 
\begin{equation}\label{equ:HC}
\text{ where } k^{HC}=\Argmax_{1\leq k\leq p q}\left \{\frac{k/p-\pi[(k)]}{\sqrt{k(p-k)}} \right \},
\end{equation}
with $\forall k =1,\dots,p$ $\pi[(k)]=2(1-\Phi(|T[(k)]|)$ and $\Phi(x)=P(\mathcal{N}(0,1)\leq x)$.  
\paragraph{Theoretical result and comments}
 
\begin{theoreme}\label{thbp}
Let $g$ be defined by Equation (\ref{gdef}) with $\hat{F}_{10}$ as given by Equation \ref{hatdeF}. 
\begin{enumerate}
\item Suppose we are using $\hat{I}^U$ as defined by Equation \ref{universal}. Assume there exists $r,R>0$ such that $0<r<\|F_{10}\|_{L_2(P_C)}\leq R$ and that $\log(p)<<\sqrt{n}$. Then, there exists $c(R)>0$ such that 
\begin{equation}\label{uppp}
\E_{P^{\otimes n}}\left [\mathcal{C}(g)-\mathcal{C}^*\right ]\leq \frac{c}{r}\left \{ \frac{\log(p)}{n}\left (1+\mathcal{R}(C^{-1/2}m_{10},n)\right )+o\left (\frac{\log(p)}{n}\right )\right \}
\end{equation}
where 
\[\mathcal{R}(C^{-1/2}m_{10},n)=\sum_{i=1}^{p}\min \left (n\frac{m^2_{10}[i]}{\sigma^2[i]},1 \right ).\]
\item Suppose now that we are using $\hat{I}^{FDR}$ as defined by Equation \ref{selectFDR}. 
Suppose that in Equation $\ref{multipletest}$ and in the definition of $\lambda_{FDR}$ below this equation, $\hat{\sigma}[i]$ equals $C[i,i]$. Define $\eta_p=p^{-\frac{1}{p}}R\sqrt{n(p)}$. If $\eta_p^q\in [\frac{\log^5(p)}{p},p^{-\delta}]$ for $\delta>0$, then, for all $0<q<2$ we have
\begin{equation}
\forall r>0,\;\;\sup_{(P_0,P_1)\in\Omega_q(R,e,r)}\E_{P^{\otimes n}}\left [\mathcal{C}(g)-\mathcal{C}^*\right ]\leq \frac{c(b_p)}{2r} \left ( \frac{\log\left (\frac{p}{R^qn(p)^{q/2}}\right )}{2Rn^{1/2}(p)}\right )^{2-q},
\end{equation}
\[c(b_p)=1+\frac{b_p}{1-b_p}+o_p(1),\]
where $b_p$ is the real value used for the choice of $k^{FDR}_{10}$, and $P^{\otimes n}$ is the law of the learning set.
\end{enumerate}
\end{theoreme}
{\it Comments about point $1$ and general comments.} The bound given by Equation \ref{uppp} can lead to a rate of convergence if one know a suitable bound for $\mathcal{R}(C^{-1/2}m_{10},n)$. These type of bounds are well known (see for example Lemma 6.1 in \cite{Candes:2006lk}) and we won't give further comment. As an example, when the number of non null components of $C^{1/2}m_{10}$ is bounded by $S$, we have 
\[\mathcal{R}(C^{-1/2}m_{10},n)\leq S,  \]
which implies that 
\[\E_{P^{\otimes n}}\left [\mathcal{C}(g)-\mathcal{C}^*\right ]=O\left (\frac{S\log(p)}{n}\right ). \]  

 The assumption $\|F_{10}\|_{L_2(P_C)}\leq R$ could be relaxed with additional technicalities in the proofs. Anyway it is easy to understand that  large values of $\|F_{10}\|_{L_2(P_C)}$ correspond to the case where the data are well separated and is not of great interest. In addition, it is often needed implicitly when one wants to bound $\mathcal{R}(C^{-1/2}m_{10},n)$. 

The assumption $\log(p)<<\sqrt{n}$ can be seen as a rather strong assumption for very large $n$. It is needed to show that the use of $\hat{\sigma}[i]$ in (\ref{universal}) gives almost the same result as the one we would have by taking $\sigma[i]$ instead.

 Note that, for certain values of $q\in]0,2[$, the rate of convergence can be fast (i.e faster than $n^{-1/2}$) under the condition that $C^{-1/2}m_{10}\in l^q$. On the other hand, assuming that $r>0$ cannot tend to zero can be seen as a margin assumption, since 
\[\|F_{10}\|_{L_2(P_C)}>r>0\implies \exists C>0\; : \; \forall \epsilon>0\; P(|1-2\eta(X)|\leq \epsilon)\leq C \epsilon. \] 
where $\eta(X)=\E[Y|X]$.
Apart from Theorem \ref{Th1} (from which Theorem \ref{thbp} can be derived) the theoretical novelty of this paper is to give upper bound on the excess risk for procedure involving a particular dimensionality reduction (Procedure $1$ and $2$ for the choice of $\hat{I}$). In Bickel and Levina \cite{bickel:2004fk} no thresholding procedure is proposed and in Fan and Fan \cite{Fans:2007vn} the choice of the threshold  is introduced after the main theoretical result to mimic the oracle bound of their Theorem $5$. In addition, most results in Fan and Fan \cite{Fans:2007vn} are established  in the case where $C=Id$. Let us recall  that  if $Y$ is a gaussian random variable with values in a Hilbert Space, then the covariance operator is necessarily nuclear. Also, the assumption used by the above mentioned authors cannot let us consider, as a limiting distribution when $p$ tends to infinity, gaussian measures with support in a Hilbert space. 

Finally, even if Theorem \ref{thbp} doesn't treat the case where $C$ is not diagonal Theorem \ref{Th1} gives hints in that direction and extending our work  with ideas from Bickel and Levina \cite{bickel:2004fk} will be the purpose of a further study.   \\

{\it Comments about point $2$.} One can use the inequality (obtained at point 4 of the comments of Theorem \ref{Th1})
\[\E[\mathcal{C}(g)-\mathcal{C}^*] \leq c \left ( \E[ \|F_{10}-\hat{F}_{10}\|^2_{L_2(P_C)}] \right )^{1/2} (\text{ for } c>0),\]

to handle the case where $\|F_{10}\|_{L_2(P_C)}$ can tend to zero when $p$ tends to infinity (no margin assumption). The rate of convergence is not anymore the fast rate.

 In point $2$, the rate of convergence is faster when $q$ is close to $0$, and slower when it is close to $2$. This leads to consider the sparsity of $C^{-1/2}(\mu_0-\mu_1)$ as a vector of $\R^p$ in a well chosen basis.

The constant $c(b_p)$ does not depend on $q\in ]0,1/2[$. We could obtain the same speed with a universal threshold ($\lambda_U=\frac{1}{n(p)}\sqrt{2\log(p)}$). In that case, the  constant $\frac{c(b_p)}{r^2}$ would be larger (cf \cite{ABDJ05}).

 In the case of the FDR reduction dimension technique the assumption about $\hat{\sigma}[i]$ is unrealistic. We do not think the result is still true without this assumption because the obtained numerical results are rather poor. Avoiding this assumption with a slight change of the procedure could be done in further work in relation with the work in \cite{ABDJ05}. 

\subsection{Numerical Results}
We present here numerical results obtained with the presented procedures. Hence, we evaluate error rate of $6$ procedures using Equation (\ref{gdef}) : 
\begin{enumerate}
\item $g^{C}$ the procedure obtained by taking $\hat{F}_{10}=\hat{C}^-(\bar{\mu}_1-\bar{\mu}_0)$ where $\hat{C}$ his the diagonal matrix with $\hat{C}[i,i]$ the empirical variance of $(X^k_j[i])_{j=1,\dots,n\; k=0,1}$. 
\item $g^U$ the procedure obtained by taking $\hat{F}_{10}$ as given by Equation \ref{hatdeF} and $\hat{I}^U$ as defined by Equation \ref{universal}.
\item $g^{FDR}$ the procedure obtained by taking $\hat{F}_{10}$ as given by Equation \ref{hatdeF} and $\hat{I}^{FDR}$ as defined by Equation \ref{selectFDR} with $q=\gamma/log(p)$ and $\gamma$ chosen by $10$-fold cross validation over an exponential grid of $\{10^0,10^{-1},\dots, 10^{-10}\}$.
\item $g^{FAIR}$ the procedure obtained by taking $\hat{F}_{10}$ as given by Equation \ref{hatdeF} and $\hat{I}^{FAIR}$ as defined by Equation \ref{equ:FAIR}.
\item $g^{Std}$ the procedure obtained by taking $\hat{F}_{10}$ as given by Equation \ref{hatdeF} and $\hat{I}^{FDR}$ as defined by Equation \ref{selectFDR} replacing the gaussian quantiles with the appropriate student quantile function, with $q=\gamma/log(p)$ and $\gamma$ chosen by $10$-fold cross validation over an exponential grid of $\{10^0,10^{-1},\dots, 10^{-10}\}$.
\item $g^{SC}$ is nearest shrunken centroid classification procedure as defined in \cite{Tibshirani:2003im}. We used the corresponding $R$ implementation in package $pamr$.
\item $g^{HC}$ the procedure obtained by taking $\hat{F}_{10}$ as given by Equation \ref{hatdeF} and $\hat{I}^{HC}$ as defined by Equation \ref{equ:HC}, with $q$  chosen by $10$-fold cross validation over a grid of $q=\{0.2,0.1,0.05,0.01\}$.
\end{enumerate}

 We made two different simulations for the numerical experiments:
\begin{itemize}
\item {\bf Simulation 1}
\[\mu_0=0,\;\;\mu_1[i]=3 \1_{i = 4 },\;\text{ and }\; C=diag(array(1,p)). \]
\item  {\bf Simulation 2}
\begin{align*}
\mu_0=0,\;\;\mu_1[1:4]=[0.01,0.5,0.02,0.5]/3,\;\mu_1[5:p]=0,\\
\;\text{ and }\; C=diag(array(c(0.01,2),p)). 
\end{align*}
\end{itemize}

where the definition of $C$ is given in R language.\\
 All the results shown in the following tables have been obtained by repeating the experiment 100 times and averaging the error rate (which are given in  \%). The corresponding $R$ code is available \footnote{robin.girard@mines-paristech.fr}. An $R$ package will be implemented in the future including more plugin type high dimensional classification procedures. 
 
 \paragraph{Simulation 1} 
 
 The results from the first experiment confirm the poor performances of the Fisher rule with respect to the other rules (which are all based on a dimensionality reduction procedure). Procedures using cross validation for tuning of the thresholding parameter perform best (procedures $g^{FDR}$, $g^{Fisher}$, $g^{HC}$ and $g^{Student}$ use cross validation). Note that standard deviation ranges from $2$ to $5$ (in the case when $n=50$ or $100$) or even $7$ to $8$ (for $n=20$). 
\begin{longtable}{|r|rrrrrr|}
 \hline
& \multicolumn{6}{c|}{$n=20$}    \\ 
 \hline
 p  & $g^{U}$ & $g^{Fisher}$ & $g^{FDR}$ & $g^{SC}$ & $g^{Std}$ & $g^{FAIR}$ \\
  \hline
100 & 13.90 & 17.27 &{\bf 7.95}  & 8.85  & 8.15   & 13.75  \\ 
 500 & 23.97& 29.82 &{\bf 8.47} & 8.75  & 8.57 & 24.85  \\ 
  \hline
 \hline
& \multicolumn{6}{c|}{$n=50$}    \\ 
 \hline
 p  & $g^{U}$ & $g^{Fisher}$ & $g^{FDR}$ & $g^{SC}$ & $g^{Std}$ & $g^{FAIR}$ \\
  \hline
100 & 8.75  & 10.39  & 6.79  & 7.08  &{\bf 6.77} & 8.19  \\ 
 500 & 14.1 & 19.87  &{\bf 7.07}  & 7.47   & 7.10  & 18.48 \\ 
  \hline
\hline
& \multicolumn{6}{c|}{$n=100$}    \\ 
 \hline
 p  & $g^{U}$ & $g^{Fisher}$ & $g^{FDR}$ & $g^{SC}$ & $g^{Std}$ & $g^{FAIR}$ \\
  \hline
100 & 7.91 & 8.89 & 6.93 & 7.05 &{\bf 6.92} & 7.55 \\ 
 500 & 10.67 & 15.31 &{\bf 7.04}& 7.05 & 7.07 & 10.44 \\ 
  \hline
\hline
\caption{Results obtained for $n=20,50,100$, $p=100,500$ with Simulation $1$. }
\label{table100}
\end{longtable}

\paragraph{Simulation $2$}
In the second simulation the signal is really hard to distinguish and there are interesting features respectively with small and large variance. The results show the importance of using cross validation. We also see that the FAIR rule (which does not use cross validation) performs better than the Universal thresholding rule especially for moderate dimension (see $n=50$ $p=100$ or $500$).  

\begin{longtable}{|r|rrrrrrr|}
 \hline
& \multicolumn{7}{c|}{$n=10$}    \\ 
 \hline
 p  & $g^{U}$ & $g^{Fisher}$ & $g^{FDR}$ & $g^{SC}$ & $g^{Std}$ & $g^{FAIR}$ & $g^{HC}$ \\
  \hline
100 & 35.6  & 37.4   & 21.85   & {\bf 20.70}   & 22.35   & 31.25 & 22.65\\
 500 & 42.5  & 44.05  & 30.6   & {\bf 25.95} & 28.85  & 38.6 & 30.95\\ 
 5000 & 47.05   & 46.9   & 39.2   & {\bf 33.75}  & 39.45  & 45.75 &  40.00\\ 
  \hline
\hline
& \multicolumn{7}{c|}{$n=20$}\\ 
 \hline
 p  & $g^{U}$ & $g^{Fisher}$ & $g^{FDR}$ & $g^{SC}$ & $g^{Std}$ & $g^{FAIR}$ & $g^{HC}$  \\
  \hline
100 &  30.22 & 32.32  &{\bf 15.17} & 15.37  & 15.32  & 23.27 & 17.80\\ 
 500 & 38.6  & 40.57 & 17.02  & {\bf 16.05} & 16.27   & 34.25  & 20.30\\ 
 5000 & 47.25  & 47.95  & 22.77   &{\bf 19.4}  & 22.77  & 45.00 & 30.45 \\ 
  \hline
\hline
& \multicolumn{7}{c|}{$n=50$} \\
 \hline
 p & $g^{U}$ & $g^{Fisher}$ & $g^{FDR}$ & $g^{SC}$ & $g^{Std}$ & $g^{FAIR}$  & $g^{HC}$  \\
  \hline
100 & 22.42  & 24.95  & {\bf12.5} & 12.91  & 12.51  & 17.15 & 16.48\\ 
 500 & 34.14  & 36.21  & 12.92 & 13.03  & {\bf12.54} & 28.66 & 16.85\\ 
 5000 & 43.69  & 45.31  & {\bf12.36} & 12.82  & 12.47  & 42.18 &  20.04\\ 
  \hline
\hline
\caption{Results obtained for $n=20,50,100$, $p=100,500,5000$ with Simulation $2$.}
\label{tableFDR}
\end{longtable}
\section{Conclusions}

We have studied the problem of discrimination in a gaussian framework of high dimension. We have shown, with finite sample lower bounds, that standard procedures fail in high dimension ($p>>n$), and have proposed procedures to resolve this problem. These procedures are based on a dimensionality reduction technique. They also can be interpreted as thresholding estimators of the normal vector $F_{10}$ to the optimal separating hyperplan : $\{x\in \R^p: \langle F_{10},x-s_{10}\rangle_{R^p}=0\}$. We have given upper bounds on the excess risk associated to these procedures that exhibit a fast rate of convergence under a sparsity assumption. These upper bounds have been derived from a general theorem (Theorem \ref{Th1}) which may bring an interest on its own for people willing to prove convergence of other procedure in the framework of linear discriminant analysis. We have provided numerical results that confirm the theoretical development of the paper. The case when $P_0$ and $P_1$ are gaussian with different covariances can be treated with similar ideas (see the author's work \cite{Girard:2008wd} but no satisfactory theoretical results exist in this case) and will be investigated in further work. The case when the covariance matrix $C$ is not diagonal will also be the purpose of a further investigation. Futur work will discuss an evaluation of robustness for the procedure with respect to non gaussian data, numerically and theoretically.
\section{Proofs}
\subsection{Fundamental Theorem}
\begin{theoreme}\label{Th1}
Suppose $g$ is given by \ref{gdef} with $\hat{s}_{10}=s_{10}$. Let us define 
\begin{equation}\label{dzero} 
d_0=\frac{1}{\|\hat{F}_{10}\|_{L_2(P_C)}}\langle \hat{F}_{10},\hat{s}_{10}-s_{10} \rangle_{\R^p}.
\end{equation}
 Then if $\alpha=\alpha_{L_2(P)}(\hat{F}_{10},F_{10})$, we have:
\begin{equation}\label{Lowerbound}
\frac{1}{2}P\left (0<\mathcal{N}(0,1)\leq \|F_{10}\|_{L_2(P_C)}\frac{1-\cos(\alpha)}{2}\right )e^{-\frac{\|F_{10}\|^2_{L_2(P_C)}}{8}}\leq \mathcal{C}(g)-\mathcal{C}^*
\end{equation}
and
\begin{equation}\label{Upperbound}
\mathcal{C}(g)-\mathcal{C}^*\leq c P\left (|\mathcal{N}(0,1)|\leq (1-\cos \alpha) \|F_{10}\|_{L_2(P_C)}\right ) +c d^2_0
\end{equation}
for a universal constant $c>0$
\end{theoreme}
{\bf Comments}
\begin{enumerate}
\item These bounds give the relation between $(\alpha,\|F_{10}\|_{L_2(P_C)}, |d_0|)$ and the excess risk. $|d_0|$ is the error term related to the estimation of $s_{10}$ and $\alpha$ is the error term related to the estimation of $F_{10}$.
\item When $\hat{s}_{10}=s_{10}$ (i.e $d_0=0$) and $\|F_{10}\|_{L_2(P_C)}$ is fixed and positive, it is necessary to have $\alpha$ tending to zero in order to have an excess risk tending to zero.  Moreover, we see that, in this case, there exists $0<C_1<C_2$ such that 
\[C_1 \alpha^2\leq \mathcal{C}(g)-\mathcal{C}^*\leq C_2 \alpha^2.\] 
\item Recall that $d=\frac{\|F_{10}\|_{L_2(P_C)}}{2}$  can be seen as a theoretical measure of the separation between $P_1$ and $P_0$ (note that the Hellinger distance can also be expressed as a function of $d$). Large values of $d$ are associated to well separated data and small values of $d$ to non separated data. Although, Inequality \ref{Lowerbound} can be used as a contribution to the problem of finding necessary condition for the separation (by a classification rule) of gaussian mixtures (such as it is treated in \cite{Chaudhuri:2008nx}).
\item If $\Pi_{F_{10}^{\bot}}$ is the orthogonal projection operator in $L_{2}(P_C)$ one can see that :
\begin{align*}
\|F_{10}\|_{L_2(P_C)}(1-\cos(\alpha))&=\|F_{10}\|_{L_2(P_C)}-\|\Pi_{F_{10}}F_{10}\|_{L_2(P_C)}\\
&\leq \min\left \{\|\Pi_{F_{10}^{\bot}}\hat{F}_{10}\|_{L_2(P_C)},\frac{\|\Pi_{F_{10}^{\bot}}\hat{F}_{10}\|^2_{L_2(P_C)}}{2\|F_{10}\|_{L_2(P_C)}} \right \}
\end{align*}
and in particular
\begin{equation}\label{expr}
\|F_{10}\|_{L_2(P_C)}(1-\cos(\alpha))\leq \frac{\|\hat{F}_{10}-F_{10}\|^2_{L_2(P_C)}}{2\|F_{10}\|_{L_2(P_C)}}.
\end{equation}
When $d_0=0$, the upper bound in this last equation is sharper than the upper bound we have by the following standard sequence of inequalities
\begin{align*}
\E[\mathcal{C}(g)-\mathcal{C}^*]&=\E[|2\eta(X)-1|\1_{g^*\neq g}] \;(\text{ with } \eta(X)=\E[X|Y])\\
 & =\E[|\psi(e^{\mathcal{L}_{10})}|\1_{g^*\neq g}]\\
 & (\text{ with } \psi(x)=\frac{1-x}{1+x}) \text{ and } \mathcal{L}_{10}=\log(\frac{dP_1}{dP_0}) ) \\
 & \leq \E[|\mathcal{L}_{10}|\1_{sign(\mathcal{L}_{10})\neq sign(\widehat{\mathcal{L}}_{10})}]\\
 & (\text{ with }\widehat{\mathcal{L}}_{10}= \frac{1}{2}\langle \hat{F}_{10},s_{10}-x\rangle_{\R^p}) \\
 & \leq \E[|\mathcal{L}_{10}-\widehat{\mathcal{L}}_{10}|]\\
 & \leq c \left ( \E[ \|F_{10}-\hat{F}_{10}\|^2_{L_2(P_C)}] \right )^{1/2} (\text{ for } c>0).
\end{align*}
which, if $\|\hat{F}_{10}\|^2_{L_2(P_C)}$ remains bounded from below (this can be seen as a margin assumption), is the square root of what can be derived from (\ref{expr}). It is also sharper than the bound given at the end of Section 2 in \cite{bickel:2004fk}.
\end{enumerate}
\subsection{Proof of Theorem \ref{thbp}}
\begin{proof} 
We separate the proof into $3$ steps.\\
 
{\bf Step 1}
First, with Theorem \ref{Th1} Equation \ref{Upperbound} and Equation \ref{expr} we have:
\[\E_{P^{\otimes n}}\left [\mathcal{C}(g)-\mathcal{C}^*\right ]\leq c \E_{P^{\otimes n}}\left [ \frac{\|\hat{F}_{10}-F_{10}\|^2_{L_2(P_C)}}{\|F_{10}\|_{L_2(P_C)}}+|d_0|^2\right ].\]

On the one hand, 
\begin{align*}
\E_{P^{\otimes_{i\in A} }}\left [ |d_0|^2\right ]&=\E_{P^{\otimes_{i\in A}}}\left [\left |\left \langle \frac{\hat{F}_{10}}{\|\hat{F}_{10}\|_{L_2(P_C)}},s_{10}-\hat{s}_{10}\right \rangle_{\R^p}\right |^2 \right ]\\
&\leq  \frac{c'}{n}
\end{align*}
for a given constant $c'$. 
On the other hand, by construction:
\[
A\hat{=}\|\hat{F}_{10}-F_{10}\|^2_{L_2(P_C)}=\|(\sigma[i]\tilde{F}_{10}[i]1_{i\in \hat{I}})_{i=1,\dots,p}-(\sigma[i]F_{10}[i])_{i=1,\dots,p}\|_{\R^p},
\]
and
\begin{align*}
\left (\sigma[i]\tilde{F}_{10}[i]\1_{i\in \hat{I}}-\sigma[i]F_{10}[i] \right )^2 & =\frac{\sigma^2[i]}{\hat{\sigma}^2[i]}\left (\frac{\bar{m}_{10}[i]}{\hat{\sigma}}\1_{i\in \hat{I}}-\frac{m_{10}[i]}{\sigma[i]}\frac{\hat{\sigma}[i]}{\sigma[i]}\right )^2\\
 &\leq  \frac{\sigma^4[i]}{\hat{\sigma}^4[i]}\left (\frac{\bar{m}_{10}[i]}{\sigma[i]}\1_{i\in \hat{I}}-\frac{m_{10}[i]}{\sigma[i]}\right )^2+ \frac{m^2_{10}[i]}{\sigma^2[i]}\left (1-\frac{\sigma[i]^2}{\hat{\sigma}[i]^2}\right )^2.
\end{align*}
Using standard inequality around the convergence of $\hat{\sigma}^2[i]$ to $\sigma^2[i]$, one can show, summing up over $i\in \{1,\dots,p\}$, that there exist a constant $c>0$ such that 
\begin{equation}\label{expA}
\E[A]\leq c\left (  \E \left[ \left\|\left ( \frac{\bar{m}_{10}[i]}{\sigma[i]}\1_{i\in \hat{I}}\right )_{i=1,\dots,p}-\left (\frac{m_{10}[i]}{\sigma[i]}\right )_{i=1,\dots,p}\right \|^2_{\R^p}\right ]+\|F_{10}\|_{L_2(P_C)}\frac{1}{n}\right ).
\end{equation}  
Hence, it only remains to bound the expectation in the right side of the preceding equation, say $\E[B]$. 
In both case (step $2$ and step $3$) we will use the fact that the covariance matrix of the vector  $C^{-1/2}\bar{m}_{10}$ equals $I_p \frac{1}{n}$.\\

{\bf Step 2 : the case of the universal procedure}\\
In the case of the universal procedure, 
\[\left ( \frac{\bar{m}_{10}[i]}{\sigma[i]}\1_{i\in \hat{I}}\right )_{i=1,\dots,p}=\left (Y[i]\1_{|\frac{Y[i]}{\sqrt{Var(Y[i])}}|\geq \frac{\hat{\sigma}[i]}{\sigma[i]}\sqrt{2\log(p)}}\right )_{i=1,\dots,p},\]
where $Y=C^{-1/2}\bar{m}_{10}$. 

Following the notations of Theorem 4 (with $n$ replaced by $p$) from Donoho and Johnstone \cite{DJ94a}, we set

\[l_p=\frac{\hat{\sigma}[i]}{\sigma[i]}\sqrt{2\log(p)}\]
$\gamma$ a positive constant, $\epsilon_p$ a positive sequence decreasing to zero
and define three different events : 
\[E=\left \{(1-\gamma)\log\log(p) \leq l_p^2-2\log(p)\leq c\epsilon_p \log(p) \right \}\]
\[ E_-=\left \{(1-\gamma)\log\log(p) \geq l_p^2-2\log(p)\right \} \;\; E_+=\left \{l_p^2-2\log(p)\geq c\epsilon_p \log(p)\right \}.\]
From the bayes formula we get: 
\[\E[B]\leq \E[B|E]+\left (P(E_+)+P(E_-)\right )\E[B].\]
We also have
\[\E[B]\leq \E[\|Y\|_{\R^p}^2]+ \|F_{10}\|^2_{L_2(P_C)}=2\|F_{10}\|_{L_2(P_C)}+2p\]
Concentration inequalities for a chi square random variables $U$ with $n-1$ degrees of freedom (see for example comments on Lemma 1 in \cite{Laurent:2000fk}) give (for $n>4$)
\begin{align*}
P(E_+)&\leq P\left ( U-(n-1)\geq \frac{n-1}{2}c\epsilon_p\right ) \\
 &\leq  P\left ( U-(n-1)\geq \sqrt{n-1} \frac{\sqrt{n-1}}{4}c\epsilon_p+ \frac{\sqrt{n-1}}{2}c \epsilon_p \right )\\
 & \leq e^{-\frac{\sqrt{n-1}}{4}c\epsilon_p}=o\left (\frac{\log(p)}{pn}\right )
\end{align*}
(because $\log(p)<<\sqrt{n}$) and 
\begin{align*}
P(E_-)&\leq P\left ( U-(n-1)\leq \frac{(\gamma-1)(n-1)}{2}c\frac{\log (p)}{\log\log(p)} \right ) \\
& \leq e^{-\frac{(\gamma-1)\sqrt{n-1}}{2}c\frac{\log (p)}{\log\log(p)}} =o\left (\frac{\log(p)}{pn}\right ).\\
\end{align*}
 This ends the proof.
  \\

{\bf Step 3 : the case of the FDR procedure}

 Theorem $1.1$ of Abramovich an .al \cite{ABDJ05},  and Theoreme $5$ point $3b.$ of Donoho and Johnstone \cite{Donoho:1994uq} then lead to the desired result.

\end{proof}
\subsection{Proof of Theorem \ref{Th1}}
In this proof, we will use the following subset of $\R^p$:
\[\hat{V}=\{x\in\R^p\;:\; \langle \hat{F}_{10},x-\hat{s}_{10} \rangle_{\R^p}\geq 0, \;\;V=\{x\in\R^p\;:\; \langle F_{10},x-s_{10} \rangle_{\R^p}\geq 0\}\]
\[\hat{V}_{2}=\{x\in\R^p\;:\; \langle C^{1/2}\hat{F}_{10},x \rangle_{\R^p}\geq d_0\;\; V_{2}=\{x\in\R^p\;:\; \langle C^{1/2}F_{10},x\rangle_{\R^p}\geq 0\}, \]
where $d_0$ is defined by Equation \ref{dzero}.
The proof is divided into four steps: in the first one we make a change of geometry and in the second one we obtain a simple expression with gaussian measure of subsets or $\R^2$. In the third one we derive the lower bound and in the fourth one the upper bound. \\

{\bf Step 1}. We have 
\begin{align*}
\mathcal{C}(g)-\mathcal{C}(g^*)&=\frac{1}{2}\left (P_0(\hat{V}\setminus V)-P_0(V\setminus \hat{V}) +P_1(V\setminus \hat{V})-P_1(\hat{V}\setminus V) \right)\\
&=\frac{1}{2}\left (P_{10}(\hat{V}\setminus V-m_{10})-P_{10}(V\setminus \hat{V}-m_{10})\right .\\
&+\left .P_{10}(V\setminus \hat{V}+m_{10})-P_{10}(\hat{V}\setminus V+m_{10})\right)
\end{align*}
where $P_{10}$ is the gaussian probability distribution with covariance $C$ and mean $s_{10}$, and $m_{10}=\frac{\mu_1-\mu_0}{2}$. Changing the geometry now gives
\begin{align*}
\mathcal{C}(g)-\mathcal{C}(g^*) = &\frac{1}{2}\left ( P(\xi-C^{1/2}F_{10}/2\in \hat{V}_2\setminus V_2)-P(\xi-C^{1/2}F_{10}/2\in V_2\setminus \hat{V}_2) \right . \\
&+ \left . P(\xi+C^{1/2}F_{10}/2\in V_2\setminus \hat{V}_2)-P(\xi+C^{1/2}F_{10}/2\in \hat{V}_2\setminus V_2) \right )
\end{align*}
where $\xi$ is a gaussian random variable on $\R^p$ with mean $0$ and covariance $I_p$. Notice that if $\alpha=\alpha_{L_2(P_{C})}(\hat{F}_{10},F_{10})$ $(\mathcal{C}(g)-\mathcal{C}(g^*))(\alpha)=(\mathcal{C}(g)-\mathcal{C}(g^*))(-\alpha)$, also, we will suppose without loss of generality that $\alpha>0$ in the rest of the proof.\\

\begin{figure}
    \center
    \includegraphics[width=9cm]{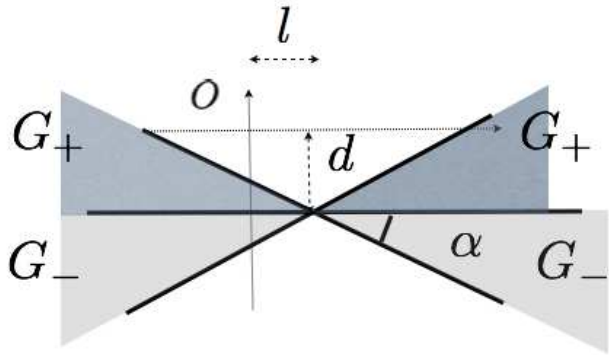}
   \caption{Figure giving the definition of  $G_+$ and $G_-$.}
   \label{MyFig1}
\end{figure}

{\bf Step 2}.
This step is roughly a geometric exercise in $\R^2$ (more precisely the span of $C^{1/2}\hat{F}_{10}$ and $C^{1/2}F_{10}$ in $\R^p$ or the span of $\langle\hat{F}_{10},.\rangle_{\R^p}$ and $\langle F_{10},.\rangle_{\R^p}$ in $L_2(P_C)$). First, it is easy to see (with step 1 result) that with a symmetry argument, we have
  \begin{equation}\label{withG}
  \mathcal{C}(g)-\mathcal{C}(g^*)=\frac{1}{2}\left (P(\mathcal{N}(0,I_2)\in G_+)-P(\mathcal{N}(0,I_2)\in G_{-})\right )
  \end{equation}
  where $G_+$ and $G_-$ are subsets of $\R^2$ defined by Figure \ref{MyFig1} with $d$ and $l$ given by:
    \[d=\frac{\|F_{10}\|_{L_{2}(P_C)}}{2} \text{ and } l=\frac{|d_0|}{\sin(\alpha)\|\hat{F}_{10}\|_{L_2(P_C)}}. \]
(note that obtaining $l$ needs a small calculation with $\R^2$ geometry).\\

\begin{figure}
    \center
    \includegraphics[width=13cm]{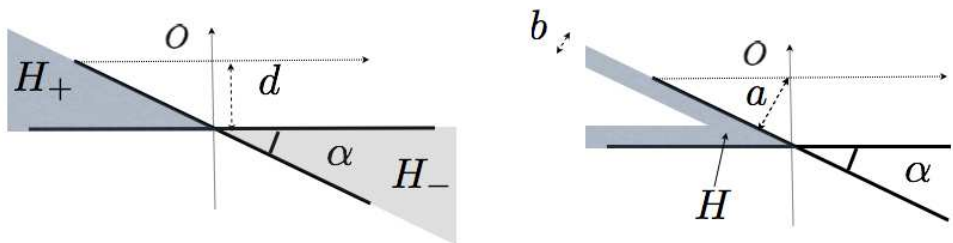}
   \caption{On the left: Figure giving the definition of  $H_+$ and $H_-$. On the Right figure defining $a$ and $b$.}
   \label{MyFig1b}
\end{figure}

{\bf Step 3 : The lower bound.}
For the lower bound, we shall first notice that 
\[\mathcal{C}(g)-\mathcal{C}(g^*)\geq \frac{1}{2}\left (P\left (\mathcal{N}(0,I_2)+\left(\begin{array}{c} l \\0\end{array}\right)\in G_+\right )-P\left (\mathcal{N}(0,I_2)+\left(\begin{array}{c} l \\0\end{array}\right)\in G_-\right )\right ),\]
and by symmetry, this gives 
\begin{equation}
\mathcal{C}(g)-\mathcal{C}(g^*) \geq P(\mathcal{N}(0,I_2)\in H_+)-P(\mathcal{N}(0,I_2)\in H_-),
\end{equation}
where $H_+$ and $H_-$ are given defined by Figure \ref{MyFig1b}.
Let $B$ be the orthogonal projection of $O$ on to the bisector of $\alpha$ in the Figure defining $H_+$ and $H_-$ (i.e Figure \ref{MyFig1b} on the left). Let us define $H=H_+\setminus S_B(H_-)$ (see Figure \ref{MyFig1b} on the right) where $S_B$ is the symmetry of center $B$ (also the symmetry of axe $(O,B)$). One can see that with this construction and the preceding equation, we have:  
\[\mathcal{C}(g)-\mathcal{C}(g^*)\geq P(\mathcal{N}(0,I_2)\in H).\]
From this equality and standard inequality on gaussian measures, we get
\[ \mathcal{C}(g)-\mathcal{C}(g^*)\geq\frac{P(\mathcal{N}(0,1)\in [0,b])}{2}e^{-\frac{(a+b)^2}{2}},\]
where $a$ and $b$ are defined by Figure \ref{MyFig1b} on the right and can be calculated easely:
\[b=(1-\cos \alpha)\frac{\|F_{10}\|_{L_2(P_C)}}{2} \;\;\;  a=\cos \alpha\frac{\|F_{10}\|_{L_2(P_C)}}{2}.\]
 This gives the announced lower bound.\\
 
 {\bf Step 4 : The upper bound.}
First, we notice that 
\[\mathcal{C}(g)-\mathcal{C}(g^*)=\mathcal{C}(g)-\mathcal{C}(\tilde{g})+\mathcal{C}(\tilde{g})-\mathcal{C}(g^*)\]
where 
\[\tilde{g}=\left\{\begin{array}{cc}1 & \text{if }\langle \hat{F}_{10},x-s_{10}\rangle_{\R^p} \geq 0\\ 0 & \text{ otherwise }\end{array}\right.\]
With step two (setting $d_0=0$), we have
\[\mathcal{C}(\tilde{g})-\mathcal{C}(g^*)\leq P(\mathcal{N}(0,I_2)\in H_+)-P(\mathcal{N}(0,I_2)\in H_-)=P(\mathcal{N}(0,I_2)\in H).\]
From this equality and standard inequality on gaussian measures, we get
\[ \mathcal{C}(\tilde{g})-\mathcal{C}(g^*)\leq P(\mathcal{N}(0,1)\in [0,2b])e^{-\frac{a^2}{2}}.\]
It now remains to bound $\mathcal{C}(g)-\mathcal{C}(\tilde{g})$. We have, following the same type of calculation we had in step 1, the following equality :  
\[\mathcal{C}(g)-\mathcal{C}(\tilde{g})=P(\xi\in [0;\epsilon])-P(\xi\in [-\epsilon;0])\]
with 
\[\xi\leadsto \mathcal{N}(m,\sigma^2),\;\;\epsilon=\sigma |d_0|,\;\; \sigma=\|\hat{F}_{10}\|_{L_2(P_C)}\text{ and }m=\langle F_{10},\hat{F}_{10}\rangle_{L_2(P_C)}.\]
Also, the desired bound follows directly from the following lemma.
\begin{Lemme}\label{lemme1}
If $m\in \R$, $\sigma>0$, $\xi\leadsto \mathcal{N}(m,\sigma^2)$ and $\epsilon>0$, then there exists $c>0$ such that
\[P(\xi\in [0;\epsilon])-P(\xi\in [-\epsilon;0])\leq c \frac{\epsilon^2}{\sigma^2}\]
\end{Lemme}
\begin{proof}

Let us call $R$ the left side of the inequality to be proved and set $\Phi(x)= P\left (\mathcal{N}(0,1)\leq x\right )$. We have 
\[R=\Phi\left (\frac{\epsilon}{\sigma}-\frac{m}{\sigma}\right )+\Phi\left (-\frac{\epsilon}{\sigma}-\frac{m}{\sigma}\right )-2\Phi\left (-\frac{m}{\sigma}\right )\]

which gives the desired result with taylor expension since there exist $C>0$ such that  $|\Phi''|\leq C$.  
\end{proof}

 \subsection{Proof of Proposition \ref{prop:un} point 1}
 \begin{proof}
 The proof is based on ideas from Bickel and Levina \cite{bickel:2004fk} used in their Theorem $1$ :
 if $C$ is the identity their exist $\xi_1,\dots,\xi_p$, $p$ $\R^p$ valued random variables forming an orthonormal basis of $\R^p$,  a random vector $(\lambda_1,\dots,\lambda_n)$ of $\R^n$ whose property are the following. 
  \begin{enumerate}
  \item The $\lambda_i$ are independent between each other, independent from $(\xi_i)_{i=1,\dots,p}$, and $n\lambda_i$ follows a $\chi^2$ distribution with $n-1$ degrees of freedom.
  \item For every $i$, $\xi_i$ is drawn in an independent and uniform fashion on the intersection of the unitary sphere of $\R^p$ and the orthogonal to $\xi_{1},\dots,\xi_{i-1}$.
  \item The empirical estimator $\hat{C}$ of $C$ verify :  
\[
  \hat{C}=\sum_{i=1}^n\lambda_i\xi_i\otimes \xi_i, 
\]
where if $x,y\in \R^p$, $x\otimes y$ is the linear operator of $\R^p$ that associate to $z\in \R^p$ the vector $\langle x,z \rangle_{\R^p}y$.  
\end{enumerate}

When $C$ not necessarily equals $I_p$, we get, $P_C-$almost-surely :
\[
C^{-1/2}\hat{C}C^{-1/2}=\sum_{i=1}^n\lambda_i\xi_i\otimes \xi_i, \text{ et }C^{1/2}\hat{C}^{-}C^{1/2}=\sum_{i=1}^n\frac{1}{\lambda_i}\xi_i\otimes \xi_i.
\]
Then, if we define $\beta_i=\langle C^{-1/2}m_{10},\xi_i\rangle_{\R^p}^2$, we have the following equations
\begin{equation}\label{bick1}
\langle F_{10},\hat{F}_{10}\rangle_{L_2(P_C)}=\langle C^{-1/2}m_{10},C^{1/2}\hat{C}^{-}C^{1/2}C^{-1/2}m_{10}\rangle_{\R^p}=\sum_{i=1}^n \frac{\beta_i}{\lambda_i},
\end{equation}
\begin{equation}\label{bick2}
\|\hat{F}_{10}\|_{L_2(P_C)}^2=\sum_{i=1}^n\frac{\beta_i}{\lambda_i^2} \text{ et }\|F_{10}\|_{L_2(P_C)}^2=\sum_{i=1}^p\beta_i.
\end{equation}
For reasons of symmetry  (the $\xi_i$ are drawn uniformly on the sphere), we have for all subset $I_n$ from $\{1,\dots,p\}$ of size $n$ :
\begin{equation}\label{bikc1}
u_{I_n,p}=\E\left [\frac{\sum_{i\in I_n}\beta_i}{\sum_{i=1}^p\beta_i}\right ]=\frac{n}{p},
\end{equation}
 From equations (\ref{bick1}) and (\ref{bick2}), if $\alpha=\alpha_{L_2(P)}(\hat{F}_{10},F_{10})$, we have ( Cauchy-Schwartz inequality ):
\[
cos(\alpha)=\frac{\sum_{i=1}^n\frac{\beta_i}{\lambda_i}}{\left (\sum_{i=1}^p\beta_i\right )^{1/2}\left (\sum_{i=1}^n\frac{\beta_i}{\lambda_i^2}\right )^{1/2}}\leq \left ( \frac{\sum_{i=1}^n\beta_i}{\sum_{i=1}^p\beta_i}\right )^{1/2}1_{|\alpha|<\pi/2}\\
\]
Hence, with Jensen inequality and Equation (\ref{bikc1}), this gives $\E[\cos(\alpha)]\leq \sqrt{\frac{n}{p}}$. This and inequality (\ref{Lowerbound}) leads to the desired result. 
\end{proof}
\subsection{Proof of proposition \ref{prop:un} point 2}

\begin{proof}
As in the preceding proposition, we are going to use Inequality (\ref{Lowerbound}). Also it is sufficient to show the following 
\[\E\left [\cos (\alpha) 1_{|\alpha|<\pi/2} \right ]\leq\frac{1}{\sqrt{p-2}}(\sqrt{n}\|F_{10}\|_{L_2(P_C)}+1).\]
We not that suffices to obtain 
\begin{equation}\label{aprouver}
\E\left [\frac{|\langle F_{10},\hat{F}_{10}\rangle_{L_2(P_C)}|}{\|F_{10}\|_{L_2(P_C)}\|\hat{F}_{10}\|_{L_2(P_C)}}\right ]\leq \frac{1}{\sqrt{p-2}}(\sqrt{n}\|F_{10}\|_{L_2(P_C)}+1).
\end{equation} 
On the other hand, 
\[
\E\left [\frac{|\langle F_{10},\hat{F}_{10}\rangle_{L_2(P_C)}|}{\|F_{10}\|_{L_2(P_C)}\|\hat{F}_{10}\|_{L_2(P_C)}}\right ] \leq \E\left [\frac{\|F_{10}\|_{L_2(P_C)}}{\|\hat{F}_{10}\|_{L_2(P_C)}}\right ]+\E\left [\frac{|\langle F_{10},\hat{F}_{10}-F_{10}\rangle_{L_2(P_C)}|}{\|F_{10}\|_{L_2(P_C)}\|\hat{F}_{10}\|_{L_2(P_C)}}\right ]
\]
\[
\hspace{1cm}\leq \E\left [\frac{\|F_{10}\|_{L_2(P_C)}^2}{\|\hat{F}_{10}\|_{L_2(P_C)}^2}\right ]^{1/2}\left (1+\E\left [\frac{\langle F_{10},\hat{F}_{10}-F_{10}\rangle_{L_2(P_C)}^2}{\|F_{10}\|_{L_2(P_C)}^2}\right ]^{1/2}\right ),
\]
where this last inequality results from Cauchy-Scwartz.
Recall that  
\[\hat{F}_{10}=F_{10}+\frac{C^{-1/2}}{\sqrt{n}}\xi,\]
where $\xi$ is a standardised gaussian random vector of $\R^p$. Also, we easily obtain,
\[\E\left [\frac{\langle F_{10},\hat{F}_{10}-F_{10}\rangle_{L_2(P_C)}^2}{\|F_{10}\|_{L_2(P_C)}^2}\right ]^{1/2}=\frac{1}{\sqrt{n}},\]
 and 
\[\frac{\|F_{10}\|_{L_2(P_C)}^2}{\|\hat{F}_{10}\|_{L_2(P_C)}^2}=\frac{\|\sqrt{n}C^{1/2}F_{10}\|^2_{\R^p}}{\|\sqrt{n}C^{1/2}F_{10}+\xi\|^2_{\R^p}}.\]
The rest of the proof follows from the following simple fact which is a consequence of Cochran Theorem and classical calculation on $\chi^2$ random variables :\\
\indent Let $\sigma>0$, $\beta\in \R^p$, $X$ a gaussian random vector of $\R^p$ with mean $\beta$ and covariance $I_p$. Then 
\[\E\left [\frac{1}{\|X\|_{\R^p}^2}\right ]\leq \frac{1}{p-2}.\] 
\end{proof}

\end{document}